\documentclass[12pt]{article}
\usepackage[latin1]{inputenc}
\usepackage{amssymb,amsmath,amsfonts}
\usepackage{color}
\newtheorem{theorem}{Theorem}
\newtheorem{proposition}[theorem]{Proposition}
\newtheorem{lemma}[theorem]{Lemma}

\newcommand{\End}{\mbox{End}}
\newcommand{\R}{\mathbb{R}}

\newcommand{\Q}{\mathbb{Q}}
\newcommand{\Sf}{\mathbb{S}}

\newcommand{\Hy}{\mathbb{H}}

\newcommand{\spa}{\mbox{span}}

\newcommand{\po}{{\hspace*{-1ex}}{\bf .  }}

\newcommand{\Lo}{\mathbb{L}}

\def\D{{\cal D}}
\def\<{{\langle}}
\def\>{{\rangle}}
\def\n{\nabla}

\def\a{\alpha}

\def\be{\begin{equation} }
\def\ee{\end{equation} }

\def\proof{\noindent{\it Proof:  }}
\def\qed{\ifhmode\unskip\nobreak\fi\ifmmode\ifinner
\else\hskip5 pt \fi\fi\hbox{\hskip5 pt \vrule width4 pt
height6 pt  depth1.5 pt \hskip 1pt }}

\begin{document}

\title{Complete minimal submanifolds with nullity\\ in the hyperbolic space}
\author{M.\ Dajczer, Th.\ Kasioumis,
A.\ Savas-Halilaj and Th.\ Vlachos}
\date{}
\maketitle

\renewcommand{\thefootnote}{\fnsymbol{footnote}} 
\footnotetext{This research is a result of the activity developed 
within the framework of the Programme in Support of Excellence Groups 
of the Regi\'on de Murcia, Spain, by Fundaci\'on S\'eneca, Science and 
Technology Agency of the Regi\'on de Murcia. This work was partially 
supported by MINECO/FEDER project reference MTM2015-65430-P and 
Fundaci\'on S\'eneca project reference 19901/GERM/15, Spain.}  
\renewcommand{\thefootnote}{\arabic{footnote}}

\begin{abstract} 
We investigate complete minimal submanifolds $f\colon M^3\to\Hy^n$ 
in  hyperbolic space with index of relative nullity at least one at 
any point. The case when the ambient space is either the Euclidean 
space or the round sphere was already studied in \cite{dksv1} and 
\cite{dksv2}, respectively. If the scalar curvature is bounded 
from below we conclude that the submanifold has to be either totally 
geodesic or a generalized cone over a complete minimal surface lying 
in an equidistant submanifold of $\Hy^n$.
\end{abstract} 

In this paper, we continue our study of a class of minimal isometric 
immersions of complete Riemannian manifolds in space forms 
$f\colon M^m\to\Q_c^n$ initiated in \cite{dksv1} for sectional curvature 
$c=0$ and continued in \cite{dksv2} for $c>0$.  The basic hypothesis 
is that the index of relative nullity satisfies $\nu\geq m-2$ 
everywhere. The goal is to conclude that under some reasonable assumption 
the submanifold has to be of a simple geometric type other than totally 
geodesic. For instance,  under the hypothesis that the Omori-Yau 
maximum principle holds on the manifold, we showed in \cite{dksv1} 
that the Euclidean submanifold has to be a $(m-2)$-cylinder.

In any of the two cases already studied, the proofs reduced to analyzing
the situation of the three dimensional  submanifolds. In fact, for 
submanifolds in spheres only  this case turned out to be possible.
This paper is devoted to complete minimal submanifolds in hyperbolic 
space $f\colon M^m\to\Hy^n$  under the  assumption that $\nu\geq m-2$.
In this ambient space this condition turns out to be quite less restrictive
than in the previously studied cases. 
Nevertheless, we have reasons to believe that the manifold being  
three-dimensional is still quite special and this is why this case 
allows a characterization of a class of submanifolds that is contained 
in the following description.\vspace{1ex}

Let $i\colon\Q_c^{n-\ell}\to\Hy^n$, $1\leq\ell\leq n-3$ and $c\neq-1$, 
denote an inclusion as a complete simply connected umbilical submanifold. 
Thus  $\Q_c^{n-\ell}$ is either a totally geodesic submanifold of 
a geodesic sphere or an equidistant  hypersurface or a horosphere 
if $c>0$, $c<0$ or $c=0$, respectively.
Let $g\colon L^2\to\Q_c^{n-\ell}$ be an isometric 
immersion of a two-dimensional Riemannian manifold.
The normal bundle of $h=i\circ g\colon L^2\to\Hy^n$ splits 
orthogonally as
$$
N_h L=i_*N_g L\oplus N_i\Q_c^{n-\ell}
$$
where $N_i\Q_c^{n-\ell} $ is regarded as a subbundle of $N_hL$.
Let $G\colon N_i\Q_c^{n-\ell}\to\Hy^n$ be defined by 
$$
G(x,w)=\exp_{g(x)}w
$$
where $\exp$ denotes the exponential map of $\Hy^n$. Then we
denote by $M^m$, $m=2+\ell$, the open subset of $N_i\Q_c^{n-\ell}$ 
where  $G$ is an immersion endowed with the metric induced by 
the map $G$.
\vspace{1ex}

\noindent{\bf Definition}:
The \emph{generalized cone} in hyperbolic space over a surface 
$g\colon L^2\to\Q_c^{n-\ell}$  is the isometric immersion
$F_g\colon M^m\to\Hy^n$, $m=2+\ell$, defined as $F_g=G|_{M^m}$.

\begin{theorem}\po\label{main}
Let $f\colon M^3\to\Hy^n$ be a minimal isometric immersion 
with index of relative nullity at least $\nu\geq 1$ at any point. 
Assume that $M^3$ is complete with scalar curvature bounded from below. 
Then $f$ is either totally geodesic or a generalized cone over a complete  
minimal surface with bounded Gauss curvature lying in an equidistant 
submanifold of $\Hy^n$.
\end{theorem}

Notice that generalized cones over minimal surfaces contained in the 
other two types  of umbilical submanifolds are not part of the theorem. 
In fact, if the surface lies inside a geodesic sphere then the generalized 
cone is never complete whereas if it lies in a horosphere then the scalar 
curvature of the cone is unbounded.

Like it happens for $c\geq 0$, in the present case where $c<0$ there 
are plenty of local examples other than generalized cones. As a matter 
of fact, a local parametrization of all minimal submanifolds 
$f\colon M^m\to\Hy^n$ with index of relative nullity $\nu=m-2$  was given 
in \cite{kas} in terms of certain elliptic space-like surfaces in either
the Lorentzian sphere or the Lorentzian flat space according to $n-m$
being even or odd, respectively. Moreover, from the results in \cite{ca},
\cite{ej} and \cite{hu} it is clear that this parametrization can be used to 
construct complete examples of any dimension other than generalized cones.

\section{Preliminaries}

Let $f\colon M^m\to\Hy^n$ denote an isometric immersion of an 
$m$-dimensional Riemannian manifold $M^m$ into hyperbolic space.
The \emph{relative nullity subspace} of $f$ at $x\in M^m$ 
is defined as 
$$
\D(x)=\{X\in T_xM\colon \alpha(X,Y)=0\;\;\text{for all}\;\;Y\in T_xM\}
$$
where $\alpha\colon TM\times TM\to N_fM$ denotes the second fundamental form
with values in the normal bundle.

The dimension $\nu(x)$ of $\D(x)$ is called the \emph{index of relative nullity} 
of $f$ at $x\in M^m$.  If $f$ has constant index of relative nullity on an open 
subset $U\subset M^m$, then $\D$ is 
integrable along $U$, its leaves are  totally geodesic submanifolds of 
$M^m$  and their images under $f$ are totally geodesics submanifolds of $\Hy^n$. 
\vspace{1,5ex}

Let $U\subset M^m$ be an open subset where the index of relative nullity 
$\nu=s>0$  is  constant. The following is a fundamental result 
in the theory of isometric immersions into space forms; cf.\ \cite{da}.

\begin{proposition}\label{comp}\po Let $\gamma\colon [0,b]\to M^m$ 
be a geodesic such that $\gamma([0,b))\subset U$ is contained in a 
leaf of relative nullity foliation. Then $\nu(\gamma(b))=s$.
\end{proposition}

We decompose any $X\in TM$ as
$$
X=X^v+X^h
$$
according to the orthogonal decomposition $TM=\D\oplus \D^{\perp}$, and
denote by $\nabla^v$ and $\nabla^h$ the components of the Levi-Civita 
connection with respect to that decomposition.

The \emph{splitting tensor} $C\colon\D\times\D^{\perp}\to\D^{\perp}$ 
is defined as
$$
C(T,X)=-{\nabla}^h_XT=-(\nabla_XT)^h
$$
for any $T\in\D$ and $X\in\D^{\perp}$. If $x\in M^n$ and $T\in\D(x)$, 
then the tensor gives rise to an endomorphism
$C_T=C(T,\cdot)\colon\D^\perp(x)\to\D^\perp(x).$
Accordingly, we regard $C$ as a map 
$$
C\colon\Gamma(\D)\to\Gamma(\End(\D^\perp)).
$$

The following differential equations are a well known consequence of 
the Codazzi equation:
\be\label{E1}
\n_S C_T=C_T\circ C_S+C_{\n^v_ST}-\<T,S\>I
\ee
where $I$ is the identity map, 
\be\label{E2}
(\nabla^h_{X}C_T)Y-(\nabla^h_{Y}C_T)X
=C_{\nabla^{v}_{X}T}Y-C_{\nabla^{v}_{Y}T}X
\ee
and 
\be\label{E3}
\nabla_T A_{\xi}|_{\D^{\perp}}=A_{\xi}|_{\D^{\perp}}\circ{C}_T
+A_{\nabla^{\perp}_T\xi}|_{\D^{\perp}}
\ee
for any $S,T\in\Gamma(\D)$, $X,Y\in\Gamma(\D^{\perp})$ and with
$$
\<A_\eta X,Y\>=\<\alpha(X,Y),\eta\>
$$ 
for any $\eta\in N_fM$. See \cite{da} or \cite{dg}.
\vspace{1ex}

A main ingredient in the proof of our result is the Omori-Yau 
maximum principle.
Recall that the Omori-Yau maximum principle is said to hold on $M^m$ 
if for any function $\varphi\in C^2(M)$ bounded from above there 
exists a sequence of points $\{x_j\}_{j\in\mathbb{N}}$ such that
$$
\varphi(x_j)>\sup\varphi-1/j,\quad \|\nabla\varphi\|(x_j)\le 1/j\quad
\text{and}\quad\Delta\varphi(x_j)\le 1/j
$$
for any $j\in\mathbb{N}$. It is well-known that the principle is 
valid if the manifold is complete and the Ricci curvature does not 
decay fast to $-\infty$. It also holds for properly immersed submanifolds 
in space forms with bounded length of the mean curvature vector field; 
cf.~\cite{amr}. 
\vspace{1ex}

The following two results that will be used for the proofs are  consequence 
of the Omori-Yau maximum principle.

\begin{proposition}\label{maxprinc2}\po
Let $M^m$ be a Riemannian manifold for which the Omori-Yau maximum principle 
holds. If $\varphi\in C^{\infty}(M)$  satisfies that $\Delta\varphi\ge 2\varphi^2$ 
then $\sup\varphi= 0$. 
\end{proposition}

\proof See \cite{amr} or \cite[Lemma 4.1]{hsv2}.\qed

\begin{proposition}\label{yau}\po
Let $M^m$ be a Riemannian manifold which Ricci curvature bounded from 
below by $-K$ for some constant $K\geq 0$. If $f\in C^\infty(M)$ is a  
harmonic function which is bounded from below on $M^m$, then
$$
\|\n f\|\leq\sqrt{(m-1)K}(f-\inf f).
$$
\end{proposition}

\proof See \cite[Theorem $3^{\prime\prime}$]{Yau}.\qed

\section{Generalized cones}

In this section, we discuss several basic properties of 
generalized cones in hyperbolic space.
\vspace{1ex}

That a smooth tangent distribution $D$ is umbilical means that 
there exists a smooth section $\delta$ of $D^\perp$ such that
$$
\<\nabla_XY,T\>=\<X,Y\>\<\delta,T\>
$$
for all $X,Y\in D$ and $T\in D^\perp$.

\begin{proposition}\label{gencones}\po Let $f\colon M^m\to\Hy^n$ be an 
isometric immersion with constant index of relative nullity  
$1\leq\nu\leq m-2$. Assume that the conullity distribution $\D^\perp$ 
is umbilical. Then $f$ is locally a generalized cone over an  isometric 
immersion $g\colon\Sigma^{m-\nu}\to\Q_c^{n-\nu}$ into an umbilical 
submanifold of $\Hy^n$. Moreover, the submanifold is globally a 
generalized cone if the relative nullity leaves are complete.
\end{proposition}

\proof Let $j\colon\Sigma^{m-\nu}\to M^m$ denote the inclusion of a leaf 
$\Sigma^{m-\nu}$ of $\D^\perp$ into $M^m$.  
Then set $h=f\circ j\colon\Sigma^{m-\nu}\to\Hy^n$. 
The normal bundle of $h$ splits as 
$$
N_h\Sigma=f_*N_j\Sigma\oplus N_fM=f_*\D|_\Sigma\oplus N_fM.
$$
By assumption, there exists $\delta\in\Gamma(\D)$ such that 
$C_T=\<T,\delta\>I$ for all $T\in \Gamma(\D)$. Thus,
$$
\tilde\nabla_X f_*T
=f_*\nabla_XT=-f_*C_TX+f_*\nabla^v_XT
=-\<T,\delta\>f_*X+f_*\nabla^v_XT
$$
for all $X\in T\Sigma$ and $T\in\Gamma(\D)$, where $\tilde \nabla$ 
is the induced connection on $f^*T\Hy^n$. Therefore the subbundle 
${\cal S}=f_*\D|_\Sigma$ of $N_h\Sigma$ is parallel with respect to the 
normal connection and the shape operator of $h$ with respect to 
any section $\eta=f_*T$ of ${\cal S}$, with $T\in\Gamma(\D)$, is given by 
$$
A^h_\eta=\<T,\delta\>I.
$$
It follows easily that there exist an umbilical inclusion 
$i\colon\Q_c^{n-\nu}\to\Hy^n$ and an isometric immersion 
$g\colon\Sigma^{m-\nu}\to\Q_c^{n-\nu}$ such that $h=i\circ g$. 
Moreover, at any $x\in\Sigma^{m-\nu}$ the fiber ${\cal S}(x)=f_*\D(x)$  
coincides with the normal space of $i$ at $i(x)$. Thus 
the generalized cone over $g$ coincides locally with $f$.
The global statement is immediate.\qed

\begin{proposition}\label{po2}\po Let $g\colon L^2\to\Q_c^{n-\nu}$
be a minimal surface. With the notation given above we have that
the following facts hold:\medskip\\
\noindent $(i)$ The generalized cone $F_g\colon M^m\to\Hy^n$, $m=2+\nu$, 
over $g$ is a minimal immersion with index of relative nullity at least $\nu$
at any point.\medskip

\noindent $(ii)$ The map $G$ is an immersion if and only if $\Q_c^{n-\nu}$ 
is a totally geodesic submanifold of either an equidistant hypersurface or 
a horosphere in $\Hy^n$. In that case $M^m$ is complete if and only if $L^2$ 
is complete. Moreover, if $\Q_c^{n-\nu}$ is contained in an equidistant 
(respectively, horosphere) hypersurface then the scalar curvature of 
$M^m$ is bounded (respectively, unbounded)  along each fiber of the 
normal bundle of  the umbilical inclusion  $i\colon\Q_c^{n-\nu}\to\Hy^n$.
\end{proposition}

\proof Let $i\colon\Q^{n-\nu}_{c}\to\Hy^n$ be a complete simply connected 
umbilical submanifold. Then let $\eta_1,\eta_2,\dots,\eta_\nu$ be a global 
orthonormal frame of the normal bundle of $i$ such that $\eta_1$ points in 
the direction of the mean curvature vector field $H$. 

Since the normal bundle $N_i\Q^{n-\nu}_c$ is a trivial vector bundle we 
have that the map $G\colon L^2\times\R^\nu\to\Hy^n$ is given parametrically 
by
$$ 
G(x,t_1,t_2,\dots,t_{\nu})=\cosh t_{\nu}\,f_{\nu-1}(x)
+\sinh t_{\nu}\,\eta_{\nu}(x),
$$
where $f_{j}$ are defined inductively by  
$f_0= g$ and 
$$
f_j=\cosh t_jf_{j-1}+\sinh t_j\eta_j,\;\;1\leq j\leq\nu.
$$
Set 
$$
h_j=\Pi_{k=j+1}^{\nu}\cosh t_k,\;\;1\leq j\leq\nu-1 
$$
and 
$$
r=h_1(\cosh t_1-\| H\|\sinh t_1).
$$
A straightforward computation gives
\begin{align*}
G_*(X)  
&=r\, g_{*}(X),\;\; X\in TL,\\
G_* (\partial_{t_j}) 
&= h_j(\sinh t_j\,f_{j-1}+\cosh t_j\,\eta_j), 
\;\; 1\le j \le \nu -1, \\ 
G_* (\partial_{t_{\nu}}) 
&=\sinh t_{\nu}\,f_{\nu-1}+\cosh t_{\nu}\,\eta_{\nu}.
\end{align*}

It is clear that the map $G$ is an immersion if and only if   
$ \|H\|\leq 1 $, which in turn is equivalent to $\Q_c^{n-\nu}$ being a 
totally geodesic submanifold of either an equidistant hypersurface or 
a horosphere in $\Hy^n$. Moreover, its second fundamental form is given by 
$$
\alpha_G(X,Y)=r\alpha_g(X,Y)
$$
if $X,Y\in TL$, and the fact that the vectors 
$\partial_{t_1},\dots,\partial_{t_\nu}$ belong to the 
relative nullity subspace. This proves part $(i)$.

The induced metric on $L^2\times\R^\nu$ is given by
$$
\<\,,\,\>_G =r^2\,\<\,,\,\>_g+\<\,,\,\>_0,
$$
where the Euclidean space $\R^\nu$ is equipped with the complete 
Riemannian metric
$$
\<\,,\,\>_0=h^2_1dt_1^2+\cdots+h^2_{\nu-1}\,dt_{\nu-1}^2+dt_{\nu}^2.
$$
It follows from Lemma 7.2 in \cite{bish} that the manifold $M^m$ is 
complete if and only if $L^2$ is complete. 

Finally, the Gauss equation yields that the scalar curvature $s$ of 
$M^m$ is given by
$$ 
s=-m(m-1)-\frac{1}{r^2}\|\alpha_g\|^2.
$$
This clearly implies that the scalar curvature  of $M^m$  is bounded 
(respectively, unbounded) along each fiber of the normal bundle of  
the umbilical inclusion  $i\colon\Q_c^{n-\nu}\to\Hy^n$ if $\Q_c^{n-\nu}$ 
is a totally geodesic submanifold of an equidistant hypersurface  
(respectively, horosphere).\qed

\section{The proof of Theorem \ref{main}}

We start the proof with a result that holds for submanifolds on any 
dimension.

\begin{lemma}\po\label{L31} Let $f\colon M^m\to\Hy^n$ be a minimal isometric
immersion with constant index of relative nullity $\nu=m-2$. Then 
$C_T\in\spa\{I,J\}$ for any $T\in\Gamma(\D)$ where $I$ is the identity 
endomorphism and $J$ denotes the almost complex structure of $\D^\perp$. 
\end{lemma}

\proof  We have from \eqref{E3} that
\be\label{eq1}
A_\xi|_{\D^{\perp}}\circ C_T=C^t_T\circ A_\xi|_{\D^{\perp}}
\ee
for any $T\in\Gamma(\D)$.
On the other hand,  the minimality condition is equivalent to 
\be\label{eq2}
A_\xi|_{\D^{\perp}}\circ J =J^t\circ A_\xi|_{\D^{\perp}}.
\ee

We first treat the hypersurface case $n=m+1$. Let $\xi$ denote a local 
unit normal vector field along $f$ and let $e_1,e_2$ be an orthonormal 
tangent frame that diagonalizes $A_\xi|_{\mathcal{D}^{\perp}}$ 
such that $J e_1=e_2$. Set
$$
u=\<\n_{e_2}e_1,T\>\;\;\text{and}\;\;v=\<\n_{e_1}e_1,T\>.
$$
From the Codazzi equations
$$
(\n_{e_1}A_\xi)e_2 =(\n_{e_2}A_\xi)e_1\;\;\text{and}\;\; (\n_{e_i} A_\xi)T
=(\n_{T}A_\xi)e_i\;\;\mbox{for}\;\; i=1,2,
$$
we obtain that
$$
u=-\<\n_{e_1}e_2,T\>\;\;\text{and}\;\;v=\<\n_{e_2}e_2,T\>.
$$
It follows that $C_T=vI-uJ$.
\vspace{1ex}

Assume now that $f$ does not reduce codimension to one. 
The normal subspaces
$$
N_1^f(x)=\spa\{\a(X,Y):\text{for all}\;X,Y\in T_xM\}
$$
have dimension at most two due to minimality.
If there is an open subset $V\subset M^m$ where $\dim N_1^f=1$, then
a simple argument using the Codazzi equation gives
that the normal subbundle $N_1^f$ is parallel in the normal connection 
along $V$. Hence $f\vert_V$ reduces codimension to one. But then due to 
real analyticity of the immersion, the same would hold globally, 
and that is a contradiction. Therefore, there is an open dense subset  
of $M^m$ where $\dim N_1^f=2$.  From (\ref{eq1}) and (\ref{eq2}) 
it now follows easily that $C_T\in\spa\{I,J\}$ as we wished.\qed
\vspace{1ex}

We now point our attention to the three-dimensional case in which 
the harmonicity established by the following result will turn out to
be fundamental.

\begin{lemma}\po\label{L8} Assume that $m=3$. Let $e_1,e_2,e_3$ be 
a local orthonormal tangent frame such that $e_3$ spans $\D^\perp$ 
and let $v,u$ are smooth functions such that $C_{e_3}=vI-uJ$. Then
\be\label{C1}
v=-\<\n_{e_1}e_3,e_1\>=-\<\n_{e_2}e_3,e_2\>\;\;\text{and}\;\;
u=\<\n_{e_1}e_3,e_2\>=-\<\n_{e_2}e_3,e_1\>.
\ee
Moreover,
\be\label{C2}
e_3(v)=v^2-u^2-1,\;\;e_3(u)=2uv,\;\;
e_1(u)=e_2(v)\;\;\text{and}\;\;e_2(u)=-e_1(v).
\ee 
Furthermore, the functions $v$ and $u$ are harmonic.
\end{lemma}

\proof The proof of \eqref{C2} is a direct 
consequence of \eqref{E1} and \eqref{E2}. The harmonicity of the 
functions $v$ and $u$ follows by a straightforward computation 
using \eqref{C1} and \eqref{C2} similar to the one given in 
\cite{dksv1} or \cite{dksv2}.\qed
\vspace{2ex}

Next we make use of the real analytic structure of a minimal 
submanifold in order to extend smoothly the relative nullity 
distribution to the totally geodesic points. 
\medskip

Let $\mathcal{A}$ denote the set of totally geodesic points of $f$. 
By Proposition \ref{comp}, the relative nullity distribution $\D$ 
is a line bundle on $M^3\smallsetminus\mathcal{A}$.
Since $f$ is real analytic we have that
$\mathcal{A}$ is a real analytic set.  According to Lojasewicz's structure 
theorem \cite[Theorem~6.3.3]{kr}, it follows that $\mathcal{A}$ locally 
decomposes as
$$
\mathcal{A}=\mathcal{V}^0\cup\mathcal{V}^1\cup \mathcal{V}^2\cup\mathcal{V}^3
$$
where each set $\mathcal{V}^k,\, 0\leq k\leq3$, is either empty or a 
disjoint finite union of $k$-dimensional real analytic subvarieties. 
A point $x_0 \in\mathcal{A}$ is called a \emph{regular point 
of dimension} $k$ if there is a neighborhood $\Omega$ of $x_0$ 
such that $\Omega\cap\mathcal{A}$ is a $k$-dimensional real analytic 
submanifold of $\Omega$. Otherwise $x_0$ is said to be a 
\emph{singular} point. Then the set of singular points is locally a 
finite union of submanifolds.
\vspace{1ex}

We can assume that $\mathcal{V}^3$ is empty since, otherwise, 
we already have by real analyticity that $f$ is a totally geodesic 
submanifold.

\begin{lemma}\po
The set $\mathcal{V}^0$ is empty.
\end{lemma}
\proof The proof goes as in \cite{dksv1} and \cite{dksv2}.\qed

\begin{lemma}\po
The set $\mathcal{V}^2$ is empty.
\end{lemma}

\proof The proof is similar with those given in \cite{dksv1} and \cite{dksv2}. 
All we have to show is that $\mathcal{V}^2$ does not contain regular points.
Suppose to the contrary and let $\Omega \subset M^3$ be an open neighborhood 
of a regular point $x_0\in\mathcal{V}^2$ such that $L^2=\Omega\cap\mathcal{A}$ 
is an embedded surface. Let $e_1,e_2,e_3,\xi_1,...,\xi_{n-3}$ be an orthonormal 
frame  adapted to $M^3$ along $\Omega$ near $x_0$. 

The Gauss map $\gamma\colon M^3\to Gr(4,n+1)$  takes values into the 
Grassmannian of oriented space-like $4$-dimensional subspaces 
in the Lorentzian space $\Lo^{n+1}$. Regarding $Gr(4,n+1)$ as a 
submanifold in $\wedge^4\Lo^{n+1}$ via the map for the Pl\"{u}cker 
embedding, we have that
$$
\gamma=f\wedge e_1\wedge e_2\wedge e_3.
$$
The coefficients of the second fundamental form are
$$
h^a_{ij}=\<\a(e_i,e_j),\xi_a\>
$$
where from now on  $1\leq i,j,k\leq 3$ and $1\leq a,b\leq n-3$.
It is easy to see that
\be\label{gm}
\gamma_*e_i=\sum_{j,a}h^a_{ij}f\wedge e_{ja}
\ee
where $e_{ja}$ is obtained by replacing $e_j$ with $\xi_a$ in 
$f\wedge e_1\wedge e_2\wedge e_3$.
Then
$$
\sum_i\<\gamma_*e_i,\gamma_*e_i\>
=\sum_{i,j,a}(h^a_{ij})^2\<f\wedge e_{ja},f\wedge e_{ja}\>=-\|\a\|^2
$$
where the inner product of two simple $4$-vectors in $\wedge^4\Lo^{n+1}$ 
is defined by
$$
\<a_1\wedge a_2\wedge a_3\wedge a_4,b_1\wedge b_2\wedge b_3\wedge b_4\>
=\det\big(\<a_i,b_j\>\big).
$$

A long but straightforward computation using the Codazzi equation yields
\be\label{laplace}
\Delta\gamma=
-\|\a\|^2\gamma+\sum_{i,a\neq b,j\neq k}h^a_{ij}h^b_{ik}\,f\wedge e_{ja,kb}
\ee
where $e_{ja,kb}$ is obtained by replacing $e_j$ with $\xi_a$ and $e_k$
with $\xi_b$ in $e_1\wedge e_2\wedge e_3$.

We identify $\wedge^4\Lo^{n+1}$ with $\Lo^N_S$ where $\binom{n+1}{4}$ 
and $S=\binom{n}{3}$ and regard $\gamma$ as a map from $M^3$ into $\Lo^N_S$. 
Denoting by $\{A_J\}_{J\in\{1,\dots,N\}}$ the corresponding base in $\Lo^N_S$, 
where $A_1,\ldots,A_S$ are timelike and the remaining vectors spacelike,
we have that
$$
\gamma=\sum_{J=1}^N w_JA_J
$$
where $w_J=-\<\gamma,A_J\>$ for $1\leq J\leq S$ and $w_J=\<\gamma,A_J\>$
for $S+1\leq J\leq N$.

We obtain from \eqref{laplace} that
\be\label{d1}
\Delta w_J= - \|\a\|^2w_J -\epsilon_J \sum_{i,a\neq b,j\neq
k}h^a_{ij}h^b_{ik}\<f\wedge e_{ja,kb},A_J\>
\ee
where
\begin{equation*}
\epsilon_J=
\begin{cases}
+ 1,\;\; 1\leq J\leq S\;\;\\ 
-1 ,\;\;  S+1\leq J\leq N.
\end{cases}
\end{equation*}

Take a local chart $\phi\colon U\to\R^3$ of coordinates $x=(x_1,x_2,x_3)$ 
on an open subset $U$ of $\Omega$ and set
\be\label{e}
e_i=\sum_j\mu_{ij}{\partial_{x_j}}.
\ee
Setting $\theta_J=w_J\circ\phi^{-1}$, we obtain the map 
$\theta\colon\phi(U)\subset\R^3\to\Lo^N_S$ 
given by
$$
\theta=\sum_{J=1}^N \theta_JA_J=(\theta_1,\dots,\theta_N).
$$
Thus $\theta=\gamma\circ\phi^{-1}$ is the representation of the 
Gauss map with respect to the above mentioned charts. From  (\ref{e}) 
and
$$
h^a_{ij}=\sum_J\<f\wedge e_{ja},A_J\>e_i(w_j)
$$
we derive that
\be\label{hh}
h^a_{ij}=\sum_{k,J}\mu_{ik}\< f\wedge e_{ja},A_J\>(\theta_J)_{x_k}.
\ee
Thus
\be\label{a}
\|\a\|^2
=\sum_{i,j,a}\Big(\sum_{k,J}\mu_{ik}\<f\wedge e_{ja},A_J\>(\theta_J)_{x_k}\Big)^2.
\ee

The Laplacian of $M^3$ is given by
$$
\Delta=\frac{1}{\sqrt{g}}\sum_{i,j}{\partial_{x_i}}
\Big(\sqrt{g}g^{ij}{\partial_{x_j}}\Big)
$$
where $g_{ij}$ are the components of the metric of $M^3$ and $g=\det(g_{ij})$. 
Using (\ref{hh}) and (\ref{a}) we see  that (\ref{d1}) is of the form
$$
\sum_{i,j}g^{ij} (\theta_J)_{x_ix_j}+C_J\big(x,\theta,\theta_{x_1},
\theta_{x_2},\theta_{x_3}\big)=0,
$$
where $C_J\colon\phi(U)\times\R^{4N}\to\R$ is given by
\begin{align*}
C_J(x,y,z_1,z_2,z_3)&=\frac{1}{\sqrt g}
\sum_{i,j}(\sqrt{g}g^{ij})_{x_i}z_{jJ}+y_{J}
\sum_{i,j,a}\Big(\sum_{k,I}\mu_{ik}\<f\wedge e_{ja},A_I\>z_{kI}\Big)^2\\
&+\epsilon_J\sum_{I,K}\sum_{\substack{i,l,m\\a\neq b,j\neq k}}
\mu_{il}\mu_{im}\<f\wedge e_{ja,kb},A_J\>\<f\wedge e_{ja},A_K\>
\<f\wedge e_{kb},A_I\>z_{mI}z_{lK}
\end{align*}
with $y=(y_1,\dots,y_N),z_i=(z_{i1},\dots,z_{iN})$, $i,m,l\in\{1,2,3\}$
and $I,J,K\in\{1,\dots,N\}$. Let  $A_{ij}=g^{ij}I_N$, $I_N$ being the 
identity $N\times N$ matrix, $C=(C_1,\dots,C_N)$ and $\bf n$ the unit 
normal field to the surface $\phi(L^2)$ in $\R^3$. Then, the vector 
valued map $\theta=(\theta_1,\dots,\theta_N)$ satisfies the elliptic
equation
$$
\mathcal{L}\theta=\sum_{i,j}A_{ij}(x)\theta_{x_ix_j}
+C\big(x,\theta,\theta_{x_1},\theta_{x_2},\theta_{x_3}\big)=0
$$
with initial conditions: $\theta$  is constant on $\phi(L^2)$ and  
${\bf n}_*(\theta)=0$  on $\phi(L^2)$.

According to the Cauchy-Kowalewsky  theorem (cf.\ \cite{t}) the above 
system has a unique solution if the surface $\phi(L^2)$ is
noncharacteristic. This latter condition is satisfied if $Q({\bf n})\neq 0$,
where $Q$ is the characteristic form given by
$$
Q(\zeta)=\det (\Lambda(\zeta))
$$
with $\zeta=(\zeta_1,\zeta_2,\zeta_3)$ and
$$
\Lambda(\zeta)=\sum_{i,j}g^{ij}\zeta_i\zeta_j I_N
$$
is the symbol of the differential operator $\mathcal{L}$. That the surface 
$\phi(L^2)$  is  noncharacteristic follows from
$$
Q(\zeta)= \Big(\sum_{i,j} g^{ij}\zeta_i\zeta_j\Big)^N.
$$
Since $C(x,y,0,0,0)=0$ the constant maps satisfy the system. 
Due to uniqueness of solutions to the Cauchy problem, we deduce that the 
Gauss map $\gamma$ is constant on an open subset of $M^3$ and that is not
possible. \qed

\begin{lemma}\po\label{f2}
The relative nullity distribution can be extended analytically over the regular 
points of the set $\mathcal{A}$.
\end{lemma}

\proof Clearly $\D$ extends continuously over the regular points of $\mathcal{A}$. 
Let $e_1,e_2,e_3$ be a local orthonormal tangent frame on an open subset $U$ 
of $M^3\smallsetminus\mathcal{A}$ as in Lemma \ref{L8}.
We view $e_3$ as a map $F\colon U\to T^1M$ into the unit tangent bundle 
of $M^3$ endowed with the Riemannian metric inherited from the Sasaki metric 
on $TM$.  We argue that the map $F=e_3$ is harmonic. In fact, from \eqref{C1} 
and \eqref{C2} we obtain that
\begin{align*}
\Delta e_3
&=\sum_{i=1}^3(\n_{e_i}\n_{e_i}e_3-\n_{\n_{e_i}e_i}e_3)\\
&=-2(u^2+v^2)e_3\\
&=-(\|\n_{e_1}e_3\|^2+\|\n_{e_2}e_3\|^2)e_3.
\end{align*}
Hence the map $F$ satisfies the differential equation
$$
\Delta F+\|\n F\|^2F =0,
$$
which is precisely the Euler-Lagrange equation for the energy 
functional of $F$ (cf.\ \cite[Proposition 1.1]{Wood}). Thus
$F\colon U\to T^1M$ is harmonic.
Since the singular set $\mathcal{A}$ has Hausdorff dimension one, 
it follows that ${\mathrm {cap}}_2(\mathcal{A})=0$.
From a result of Meier \cite[Theorem $1$]{me} it follows that $F$ 
is of class $C^2$ on $U$. But then $F$ is real analytic by a 
result due to Eells  and Sampson \cite[Proposition  p.\ 117]{ee}.\qed

\begin{lemma}\po
The set $\mathcal{A}$ has no singular points.
\end{lemma}

\proof It now follows immediately using Proposition \ref{comp}.
\vspace{2ex}\qed

\noindent \emph{Proof of Theorem \ref{main}:}
We have seen that the relative nullity distribution $\D$ extends 
to a global line bundle, also denoted by $\D$. By passing to the 
$2$-fold covering, if necessary, we have that this line bundle is 
trivial. Thus it is spanned by a globally defined unit section $e$. 
Hence, there is a unique, up to sign, orthogonal almost complex 
structure $J\colon\D^\perp\to\D^\perp$. By Lemmas \ref{L31} and 
\ref{L8} there are harmonic functions $u,v\in C^{\infty}(M)$ such that 
$$
C_e=vI-uJ.
$$

To obtain the proof of the theorem all we have to show is that $u$
vanishes. In fact, if that is the case then the result will follow from 
Propositions \ref{gencones} and \ref{po2}.

Making use of the first two equations in \eqref{C2} and that the functions
$u,v$  are harmonic, we obtain that
\begin{align*}
\Delta(u^2+v^2-1)
&=2\|\nabla u\|^2+2\|\nabla v\|^2\\
&\geq 2(e(u))^2+2(e (v))^2\\
&= 8u^2v^2+2(v^2-u^2-1)^2\\
&\geq 2(u^2+v^2-1)^2.
\end{align*}
Since the Ricci curvature of $M^3$ is bounded from below,
then Proposition \ref{maxprinc2} applies and gives that
$u^2 + v^2\le 1$. Hence $u$ and $v$ are bounded functions.

We claim that $v^2<1$. Suppose to the contrary that there is $x_0\in M^3$ 
such that $|v(x_0)|=1$. The maximum principle for harmonic functions yields that 
$v=1$ or $v=-1$ everywhere. Hence $C_e=\pm I$.
We have using \eqref{E3} that
$$
e(\|\alpha\|^2)=e\big(\sum_{j=1}^{n-3}\mbox{tr}(A^2_{\xi_j})\big)  
=\sum_{j=1}^{n-3}\mbox{tr}(\nabla_e A^2_{\xi_j}) 
=2\sum_{j=1}^{n-3}\mbox{tr}(A_{\xi_j}\circ C_e\circ A_{\xi_j})
=\pm 2\|\alpha\|^2
$$
where $\xi_1,\dots,\xi_{n-3}$ is an orthonormal normal frame parallel along 
$\gamma$. Thus
$$
\|\alpha(\gamma(t))\|^2=ce^{\pm t}
$$
where $c>0$ is a constant. Therefore $\|\alpha\|$ is unbounded along 
$\gamma$. This clearly contradicts the assumption on the scalar 
curvature and proves the claim. 

Let $\gamma\colon\R\to M^3$ be a unit speed geodesic contained in 
a leaf of the relative nullity foliation. 
Since $v^2<1$, we have from the first equation in \eqref{C2} that
$$
(v\circ\gamma)'=(v\circ\gamma)^2-(u\circ\gamma)^2-1\leq (v\circ\gamma)^2-1.
$$
Hence the function $v\circ\gamma\colon\R\to(-1,1)$ is strictly 
decreasing and satisfies $\sup v\circ\gamma=1$ and $\inf v\circ\gamma=-1$.  
Thus the function $v$ changes sign only once along each leaf of the relative 
nullity foliation. From the first equation in \eqref{C2} and  $v^2<1$ 
it follows that
$$
e(v)=v^2-u^2-1<0.
$$
Since $0$ is a regular value of $v$, the level set $L^2=v^{-1}(0)$ is 
a $2$-dimensional connected submanifold of $M^3$ and 
the map $\rho\colon L^2\times\R\to M^3$ defined by 
$$
\rho(x,t)=\mbox{exp}_xte(x) 
$$
is a diffeomorphism. Consider the smooth function $\phi\colon L^2\times\R\to\R$ 
given by
$$
\phi\circ\rho^{-1}=
\dfrac{-2v}{1+u^2+v^2+\sqrt{(1+u^2+v^2)^2-4v^2}}\cdot
$$
Setting $\psi=\phi\circ\rho^{-1}$, we have that
\be\label{46}
\frac{\psi}{1+\psi^2}=\frac{-v}{1+u^2+v^2}\cdot
\ee
A straightforward computation using \eqref{C2} yields
\be\label{45}
e(\psi)=1-\psi^2.
\ee
Since $\phi$ vanishes on $L^2$ we obtain that
$\phi(x,t)=\tanh t$. Thus $\psi$ is bounded on $M^3$.
Hence $\theta\in C^\infty(M)$ given by
$$
\theta=u^2+(v+\psi)^2
$$
is also bounded. 
Using \eqref{C2} and \eqref{45} we readily see that 
\begin{align}\label{50}
e(\theta)
&=2u e(u)+2(v+\psi)(e(v)+e(\psi))\nonumber\\
&=4u^2v+2(v+\psi)(v^2-u^2-\psi^2)\nonumber\\
&=2(v-\psi)\theta.
\end{align}
Since  $u$ and $v$ are harmonic functions, we obtain that
\begin{align}\label{48}
\Delta \theta 
&=2\|\nabla u\|^2+2(v+\psi)
\Delta\psi+2\|\nabla(v+\psi)\|^2\notag \\
&\geq  8u^2v^2+2(v+\psi)\Delta\psi+2(e(v)+e(\psi))^2  \notag \\
&=8u^2v^2+2(v+\psi)\Delta\psi+2(v^2-u^2-\psi^2)^2.
\end{align}
On the other hand, it follows from \eqref{46} that  
\be\label{49}
\frac{(1-\psi^2)(1+u^2+v^2)^2}{(1+\psi^2)^2}\nabla\psi 
=2uv\nabla u-(1+u^2-v^2)\nabla v.
\ee
Using  the harmonicity of $u$ and $v$ again, a straightforward 
computation  gives 
\begin{align}\label{410}
\frac{1-\psi^2}{2(1+\psi^2)^2}\Delta\psi 
=&\,\frac{v(1-3u^2+v^2)}{(1+u^2+v^2)^3}
\|\nabla u\|^2+\frac{2u(1+u^2-3v^2)}
{(1+u^2+v^2)^3}\<\nabla u,\nabla v\>\nonumber\\
&+\frac{v(3+3u^2-v^2)}{(1+u^2+v^2)^3}
\|\nabla v\|^2+\frac{3\psi-\psi^3}{(1+\psi)^3}
\|\nabla\psi\|^2.
\end{align}

Since $\theta$ is bounded, by the Omori-Yau maximum principle 
there is a sequence $\{x_{j}\}_{j\in\mathbb{N}} $ of points in $M^3$ 
such that 
\be\label{omori}
(i)\;\lim \theta(x_j)=\sup \theta,\;\;(ii)\;\|\nabla \theta(x_j)\|\leq 1/j\;\;
\text{and}\;\;(iii)\;\Delta \theta(x_j)\leq 1/j.
\ee
Taking a subsequence, we have that $\lim u(x_j)=u_0,\lim v(x_j)=v_0$
and $\lim\psi(x_j)=\psi_0$.
Estimating at $x_j$ and letting $j\rightarrow \infty$, we obtain from $(i)$
and $(ii)$ of \eqref{omori} and \eqref{50} that
$$
(v_0-\psi_0)\sup \theta=0.
$$
We conclude that $u$ has to vanish unless $v_{0}=\psi_0$.

Suppose now that $v_{0}=\psi_0$. We have from \eqref{46} that $v_0=\psi_0=0$.
On the other hand, since the Ricci curvature of $M^3$ is bounded from below 
it follows from Proposition \ref{yau}
that $\|\nabla u\|$ and $\|\nabla v\|$ are bounded. Hence, from
\eqref{49}, \eqref{410} and since $\psi_0=0$, we have that 
$\Delta\psi(x_j)$ is bounded. Passing to the limit and using part
$(iii)$ of \eqref{omori}, we obtain from \eqref{48} that $u_0=0$. 
It follows using  part $(i)$ of \eqref{omori} that
$\sup \theta=0$.
Thus the function $u$ vanishes, and this concludes the proof.\qed

\noindent Marcos Dajczer\\
IMPA -- Estrada Dona Castorina, 110\\
22460--320, Rio de Janeiro -- Brazil\\
e-mail: marcos@impa.br
\medskip

\noindent Theodoros Kasioumis\\
University of Ioannina \\
Department of Mathematics\\
Ioannina--Greece\\
e-mail: theokasio@gmail.com
\medskip

\noindent Andreas Savas-Halilaj\\
University of Ioannina \\
Department of Mathematics\\
Ioannina--Greece\\
e-mail: savasha@yahoo.com
\medskip

\noindent Theodoros Vlachos\\
University of Ioannina \\
Department of Mathematics\\
Ioannina--Greece\\
e-mail: tvlachos@uoi.gr
\end{document}